# Comparaison de deux outils pour les analyses de disponibilité de production : MAROS et GRIF/BStok

# Comparison of two tools for production availability analyses: MAROS and GRIF/BStok

P. Buvry, F. Brissaud, B. Declerck, et H. Varela
DNV France
69 rue Chevaleret, 75013 Paris, France
Tel. : +33 (0)1 83 94 05 01
E-mail : pierrick.buvry@dnv.com

**Résumé**
Cette comparaison a pour objet de recenser différentes fonctionnalités offertes par deux logiciels parmi les plus utilisés pour les analyses de disponibilité de production : MAROS (développé par DNV) et GRIF/BStok (développé par Total). Les différentes fonctionnalités sont présentées et comparées sous forme de tableaux. Un cas d'étude permet ensuite de comparer les résultats obtenus par l'utilisation des deux outils. MAROS et GRIF/BStok offrent tous deux des analyses de disponibilité de production qui répondent aux besoins industriels. Si l'utilisation de MAROS offre plus de « convivialité », une contrepartie est l'aspect « boite noire » qui peut en résulter. BStok, grâce à l'exploitation sous-jacente des réseaux de Petri, offre une plus grande flexibilité de modélisation, moyennant la nécessité de disposer d'un plus haut degré d'expertise.

**Summary**
This comparison aims to present different functionalities provided by two software tools among the most used for production availability analyses: MAROS (developed by DNV) and GRIF/BStok (developed by Total). The different functionalities are presented and compared through tables. Then, a case study allows comparing the results obtained by the two tools. MAROS and GRIF/BStok both provide production availability analyses that meet industrial needs. While the use of MAROS brings to the user more "conviviality", a drawback is the "black box" aspect that results from it. BStok, thanks to the underlying exploitation of Petri nets, supplies more modeling flexibility in exchange of higher expertise.

## I. Introduction

La viabilité des exploitations (pétrolières, gazières, chimiques) dépend en grande partie de la disponibilité des systèmes et procédés mis en jeux. La *disponibilité* est définie normativement par *l'aptitude d'une entité à être en état d'accomplir une fonction requise dans des conditions données, à un instant donné ou* [en moyenne] *pendant un intervalle de temps donné, en supposant que la fourniture des moyens nécessaires est assurée* (CEI 2001). Cette définition est basée sur le temps, et sur un état unique (i.e. l'état « opérant ») d'une entité. Parce que les systèmes de production peuvent être dans une multitude d'états évoluant continûment entre « production nulle » et « production maximale », l'état « opérant » doit donc être défini en fonction d'un niveau de production de référence (e.g. niveau de capacité nominale, production contractuelle).

Cependant, parce que cette disponibilité ne différencie pas les états où la production est légèrement ou fortement éloignée d'un niveau de référence et ce, même si les effets sur la production globale peuvent être importants, cette définition est trop restrictive pour l'évaluation des performances des systèmes de production. D'autres mesures de disponibilité ont alors été proposées, notamment celles présentées par T. Aven en 1987 (Aven, 1987). La mesure retenue par la norme ISO 20815 (ISO, 2008) (qui dérive de la norme norvégienne Norsok Z-016 (NOSROK, 1998)) est la *disponibilité de production*, définie par *le rapport entre la production et la production planifiée, ou tout autre niveau de référence, sur une période de temps spécifiée*. (Cette définition est ainsi basée sur un volume et non un état, et la mesure qui en résulte n'est donc pas une probabilité.)

Les analyses de disponibilité de production sont, en particulier, utilisées pour :
- prédire les performances de production, et vérifier le respect des objectifs et des exigences ;
- identifier les conditions opérationnelles, sous-systèmes et équipements les plus critiques, et trouver des mesures d'amélioration des performances ;
- comparer différentes options et permettre la sélection / optimisation d'équipements, de configurations, et d'actions de maintenance en y intégrant des considérations économiques.

Face à la complexité des analyses de disponibilité de production, les modèles « comportementaux » couplés à des simulations de Monte Carlo (Signoret *et al.*, 2004 ; Boiteau *et al.*, 2006 ; Chang *et al.*, 2010) sont souvent privilégiés aux modèles « à états » (Aven T., 1985 ; Zio E. *et al.*, 2006). C'est notamment le cas pour deux des outils parmi les plus utilisés pour ce type d'analyses :
- MAROS (« Maintainability Availability Reliability Operability Simulation »), développé par DNV (DNV, 2012) ;
- BStok (module de la suite GRIF), développé par TOTAL (GRIF, 2012).

Ces deux logiciels utilisent des blocs diagrammes (stochastiques) de flux afin de modéliser un système de production, et des simulations de Monte Carlo pour effectuer les analyses. MAROS est directement basé sur un générateur aléatoire d'évènements (« event-driven »), tandis que BStok exploitent des réseaux de Petri et le moteur de calcul MOCA-RP. En pratique, MAROS et BStok possèdent tous deux des caractéristiques bien particulières. La comparaison de ces caractéristiques est l'objet de cette communication.

La Section II présente brièvement MAROS et de BStok. La Section III propose une comparaison de certaines fonctionnalités offertes par ces deux outils. La Section IV introduit un cas d'étude afin d'illustrer la mise en application de MAROS et de BStok, et de comparer les résultats obtenus. Enfin, la Section V conclut cette communication.

## II. Présentation de MAROS et de BStok

### II.1. Présentation de MAROS

MAROS est l'acronyme de « Maintainability, Availability, Reliability, Operability Simulation ». Les premiers prototypes ont été développés par la société britannique Jardine Technology Ltd. en 1983. En 2008, cette société a été rachetée par Det Norske Veritas (DNV), et c'est maintenant DNV Software qui assure le développement et la distribution du logiciel.

MAROS fut à l'origine développé pour l'évaluation des performances des unités « amont » (i.e. traitement des « produits bruts ») de production de pétrole et de gaz (e.g. unités offshore, systèmes sous-marins, unités flottantes). Les utilisateurs de ce logiciel se situent ainsi dans tous les principaux pays producteurs et exploitants dans le domaine du pétrole et du gaz. Cependant, MAROS est aujourd'hui aussi utilisé dans d'autres secteurs comme la pétrochimie et l'industrie des procédés.

Parmi les principales caractéristiques de MAROS, on trouve : une interface graphique assez travaillée ; une définition de la capacité des puits (i.e. sources des flux) et de la demande de production (e.g. production contractuelle) au cours du temps ; une gestion multi-produit (e.g. pétrole, gaz, eau) ; des mécanismes de compensation (e.g. stockage) ; la prise en compte de nombreuses caractéristiques de fiabilité des équipements, et de maintenance ; la définition d'évènements logiques (i.e. évènements conditionnels) ; l'inclusion de caractéristiques propres aux activités de production de pétrole et de gaz (e.g. envoi de gaz à la torche ; « rampes de production ») ; certaines analyses de coûts ; et des analyses de production relativement complètes (i.e. disponibilités, contributeurs à l'indisponibilité) fournies avec une interface graphique bien détaillée.

La modélisation sous MAROS est basée sur des blocs diagrammes (stochastiques) de flux, et les analyses sont directement réalisées par des simulations de Monte Carlo (avec un générateur aléatoire d'évènements, dit « event-driven »). Enfin, un certain nombre de « patchs » est disponible, permettant de modifier certaines caractéristiques du logiciel. Dans la présente communication, nous nous limiterons exclusivement aux fonctionnalités « de base » (i.e. sans les « patchs ») de MAROS.

### II.2. Présentation de BStok

BStok est un module de la suite GRIF, une plate-forme logicielle développée par Total. GRIF est le fruit de plus 25 ans de recherche et de développement en sûreté de fonctionnement. Aujourd'hui, le développement et la distribution du logiciel sont assurés via la société SATODEV.

La construction d'un modèle avec BStok se fait par l'intermédiaire de « prototypes » dont les plus basiques permettent de modéliser : des équipements avec de nombreuses caractéristiques (e.g. capacités de production, divers modes de défaillance, états passifs) ; des équipes de réparation (avec diverses contraintes) ; et des pièces de rechange (avec une certaine logistique). Des caractéristiques plus avancées sont disponibles grâce à l'utilisation de bibliothèques spécifiques. Chaque élément forme un bloc, et le routage des flux se fait par des connecteurs et des liens entre ces blocs.

L'utilisation de chaque « prototype » génère automatiquement la création d'un réseau de Petri (stochastique et à prédicats) qui est masquée « derrière » la modélisation en blocs diagrammes (stochastiques) de flux. Les analyses sont alors effectuées par le simulateur MOCA-RP, reconnu pour ses performances d'analyses avec les réseaux de Petri. Un autre avantage de cette approche et que des « passerelles » sont possibles avec les réseaux de Petri sous-jacents afin de rendre le modèle entièrement transparent et modifiable par les initiés. Dans la présente communication, nous nous limiterons exclusivement aux fonctionnalités « de base » (prototypes par défaut) de BStok.

## III. Comparaison des fonctionnalités de MAROS et de BStok

### III.1. Modélisation des évènements

**Table 1.** Évènements

|  | MAROS | BStok |
|---|---|---|
| Défaillances des équipements de production | - modes de défaillance illimités<br>- modes de défaillance indépendants<br>- réglage d'une défaillance à la sollicitation pour n'importe quel élément et mode de défaillance<br>- modélisation d'évènements conditionnels (« simultanés », « concurrent », « à la sollicitation », etc.), limités aux éléments actifs | - deux modes de défaillance : dégradé et critique<br>- réglage d'un taux de transition depuis le mode dégradé vers le mode critique<br>- réglage d'une défaillance à la sollicitation, limitée aux éléments passifs et aux modes critiques<br>- modélisation de modes communs de défaillance, limitées aux modes critiques |
| Défaillances des équipements de sécurité | - pas de modélisation spécifique liée aux équipements de sécurité (toutes les défaillances sont détectées en ligne) | - modélisation spécifique aux équipements de sécurité (modes de défaillances cachées et/ou révélées) |
| Défaillances systématiques | - activations/inhibitions possibles d'évènements/ défaillances à des instants programmés | - définition de l'état initial de chaque composant (marche, mode dégradé, ou mode critique) |
| Arrêts planifiés | - programmation d'arrêts planifiés selon des fréquences et/ou des dates<br>- réglage de la durée des arrêts planifiés<br>- synchronisation possible d'évènements avec des arrêts planifiés | - programmation d'arrêts planifiés selon des fréquences<br>- réglage de la durée des arrêts planifiés<br>- pas de synchronisation particulière avec des arrêts planifiés |
| Évènements climatiques | - modification des taux de réparation (mais pas des taux de défaillance) en fonction des saisons | - pas de modélisation spécifique liée aux évènements climatiques |

## III.2. Modélisation des équipements

**Table 2.** Équipements

|  | MAROS | BStok |
|---|---|---|
| **Capacités de production** | - réglage de la capacité nominale des composants<br>- réglage de la capacité perdue lors de la défaillance<br>- réglage de la capacité perdue lors de la réparation | - réglage de la capacité nominale des composants<br>- réglage de la capacité perdue lors de la défaillance, limitée aux modes dégradés<br>- capacité perdue lors de la défaillance en mode critique fixée à 100%<br>- capacité perdue lors de la réparation fixée à 100% |
| **Taux de transition** | - réglage des taux de défaillance et de réparation parmi de nombreuses lois<br>- pas de réglage de dépendance entre les taux de défaillance/réparation et les flux entrant, ainsi que les capacités des autres éléments | - réglage des taux de défaillance et de réparation parmi de nombreuses de lois<br>- réglage de dépendances entre les taux de défaillance/réparation et les flux entrant, ainsi que les capacités des autres éléments |
| **Rampes de production** | - modélisation spécifique liée aux rampes de production (i.e. durée requise pour passer d'une capacité nulle à une capacité maximale) | - pas de modélisation spécifique liée aux rampes de production |

## III.3. Modélisation des flux / configurations

**Table 3.** Flux / configurations

|  | MAROS | BStok |
|---|---|---|
| **Flux de produits** | - nombre de sources illimité<br>- gestion jusqu'à trois produits, avec différentiation des produits limitée aux « blocs parallèles » | - nombre de sources illimité<br>- gestion d'un seul produit |
| **Redondance** | - création de redondances actives ou passives<br>- réglage de la capacité de chaque branche des « blocs parallèles » | - création de redondances actives ou passives<br>- réglage de la capacité de chaque branche des « blocs parallèles » |
| **Stockage** | - possibilité limitée de modélisation de « buffer » | - pas de modélisation spécifique liée au stockage |
| **Envois à la torche** | - réglage de contraintes d'envoi à la torche en fonction du temps et/ou de la quantité envoyée<br>- réglage des ratios de produits envoyés à la torche | - pas de modélisation spécifique liée à la torche |
| **Changements de configurations** | - activations/inhibitions possibles de branches à des instants programmés<br>- changements de configurations possibles des « blocs parallèles » à des instants programmés | - pas de modélisation spécifique liée aux changements de configurations |

## III.4. Modélisation de la maintenance

**Table 4.** Maintenance

|  | MAROS | BStok |
|---|---|---|
| **Équipes de maintenance** | - réglage du nombre d'équipes de maintenance et de réparateurs par équipe<br>- réglage des attributions des équipes<br>- réglage de plages horaires de disponibilité, de temps de mobilisation, de nombres de jours de disponibilité, d'emplacements des équipes<br>- définition de coûts d'équipes de maintenance | - réglage du nombre d'équipes de maintenance et de réparateurs par équipe<br>- réglage des attributions des équipes<br>- réglage de plages horaires de disponibilité et de temps de mobilisation<br>- pas de modélisation spécifique liée aux coûts d'équipes de maintenance |
| **Ressources de maintenance** | - définition de diverses ressources (e.g. pièces de rechange, outils)<br>- réapprovisionnement des ressources de maintenance à la demande ou fréquentiel<br>- réglage de stocks initiaux, de seuils de réapprovisionnement, de temps de réapprovisionnement, de temps de mobilisation, et d'emplacements de ressources de maintenance<br>- définition de coûts de ressources de maintenance | - définition de pièces de rechange<br>- réapprovisionnement des ressources de maintenance à la demande ou fréquentiel<br>- réglage de stocks initiaux, de seuils de réapprovisionnement, de temps de réapprovisionnement<br>- pas de modélisation spécifique liée aux coûts de ressources de maintenance |
| **Maintenances préventives** | - programmation de maintenances préventives selon des fréquences et/ou des dates, limitées aux éléments actifs<br>- réglage de la capacité perdue lors des maintenances préventives | - programmation de maintenances préventives selon des fréquences<br>- capacité perdue lors des maintenances préventives fixée à 100% |
| **Priorités de maintenance** | - trois types de priorité de maintenance, avec différentes caractéristiques<br>- pas de modélisation spécifique liée aux conditions de réparation | - pas de modélisation spécifique liée aux priorités de maintenance<br>- réglage de conditions de réparation, limitées aux modes dégradés |

# IV. Application de MAROS et de BStok

## IV.1. Description du cas d'étude

Le cas d'étude présenté ici porte sur un extrait de système de production de gaz, constitué des parties de séparation, de compression, et de traitement. Il inclut aussi des systèmes de refroidissements et des vannes de sécurité (ESDV). Ce système a volontairement été simplifié, tout en cherchant à y inclure de nombreux points d'intérêts au regard des fonctionnalités présentées dans la Section III (e.g. redondances passives et actives, causes communes de défaillance, défaillances à la sollicitation, maintenances préventives et correctives, équipes de maintenance, pièces de rechange, arrêts planifiés). Les caractéristiques du système ont été sélectionnées afin de ne pas privilégier les fonctionnalités d'un outil plus que l'autre, ce qui a nécessité un certain nombre de simplifications (e.g. maximum de deux modes de défaillance par équipement, défaillances à la sollicitation limitées aux éléments passifs et aux modes critiques, causes communes de défaillance limitées aux modes critiques, détection de toutes les défaillances en ligne, pertes de capacité de 100% lors des défaillances en mode critique et lors des réparations, pas de « rampe de production », pas de stockage ni de contraintes d'envoi à la torche, maintenances préventives limitées aux éléments actifs).

La description des équipements et des actions de maintenance figure respectivement dans les Tables 5 et 6.

**Table 5.** Description des équipements

| | Hypothèse sur les composants | | Perte de capacité de production | | Taux (heure$^{-1}$) | |
|---|---|---|---|---|---|---|
| Equipement | Architecture | Modes de défaillance | Avant réparation | Pendant réparation | Taux de défaillance | Taux de réparation |
| ESDV | 4 en série sur ligne principale | Critique | 100% | 100% | 10,8×10$^{-6}$ | 1/2 |
| Séparateur | 1 en série sur ligne principale | Dégradé | 50% | 100% | 409,6×10$^{-6}$ | 1/5 |
| | | Critique | 100% | 100% | 143,6×10$^{-6}$ | 1/5 |
| Compresseur (Bloc de 3 branches en parallèle) | 3*50% dont 1 passif sur ligne principale | Dégradé | 50% | 100% | 1 084,4×10$^{-6}$ | 1/10 |
| | | Critique | 100% | 100% | 1 080,8×10$^{-6}$ | 1/17 |
| | Défaillance de cause commune (DCC) de 5% pour le mode critique Défaillance critique à la sollicitation avec une probabilité de 0,05, pour l'élément passif | | | | | |
| Système de refroidissement | 2*100% sur chacune des 3 branches en parallèle | Dégradé | 50% | 100% | 657,6×10$^{-6}$ | 1/4 |
| | | Critique | 100% | 100% | 82,4×10$^{-6}$ | 1/4 |
| | Défaillance de cause commune (DCC) de 5% pour le mode critique, par couple d'éléments sur chacune des branches parallèles actives | | | | | |
| Système de traitement | 2*100% dont 1 passif sur ligne principale | Dégradé | 50% | 100% | 119,2×10$^{-6}$ | 1/13 |
| | | Critique | 100% | 100% | 365,6×10$^{-6}$ | 1/46 |
| | Défaillance de cause commune (DCC) de 5% pour le mode critique Défaillance critique à la sollicitation avec une probabilité de 0,05, pour l'élément passif | | | | | |

**Table 6.** Description de la maintenance

| Équipements | Équipe affectée (correctif et préventif) | Maintenance préventive avec perte de capacité de production de 100% | |
|---|---|---|---|
| | | Fréquence | Durée |
| ESDV | Équipe A | Tous les 4 ans | 2 heures |
| Séparateur | Équipe A | Tous les 4 ans | 4 heures |
| Compresseurs | Équipe B | Tous les 2 ans | 12 heures |
| Système de refroidissement | Équipe A | Tous les 4 ans | 2 heures |
| Système de traitement | Équipe A | Tous les 2 ans | 12 heures |

Hypothèses additionnelles :
- La durée de vie du système est de 20 ans (il s'agit de la durée sur laquelle les analyses de disponibilité de production sont réalisées).
- Un arrêt général d'une durée de 10 jours est planifié tous les 4 ans (« arrêt quadriennal »), avec une perte de capacité de production globale de 100%.
- Les maintenances préventives ne sont modélisées que pour les éléments actifs (les maintenances préventives des éléments passifs sont supposées ne pas entraîner de pertes de production).
- Les maintenances préventives des compresseurs ont lieu séquentiellement (elles nécessitent la mobilisation d'une même équipe). Il en est de même pour l'ensemble des autres équipements.
- Lorsqu'elles ont lieux la même année, les maintenances préventives sont réalisées pendant l'arrêt général.
- Un stock d'une pièce de rechange par équipement est considéré (hors ESDV), avec un temps de réapprovisionnement à la demande de 72 heures.
- Toutes les réparations entrainent la consommation d'une pièce de rechange, sauf pour les ESDV.
- L'Équipe A de maintenance a un temps de mobilisation de 1 heure, et l'Équipe B de maintenance a un temps de mobilisation de 24 heures.

## IV.2. Application de MAROS au cas d'étude

La Figure 1 présente le bloc diagramme de flux du modèle MAROS. Chaque branche du bloc parallèle « compresseurs/refroidissement » est constituée d'un compresseur en série avec un sous-bloc parallèle composé de deux unités de refroidissement en redondance active, tel que présenté sur la Figure 2. Pour chaque équipement, les modes de défaillance sont ensuite définis tel que présenté sur la Figure 3 (exemple d'un compresseur).

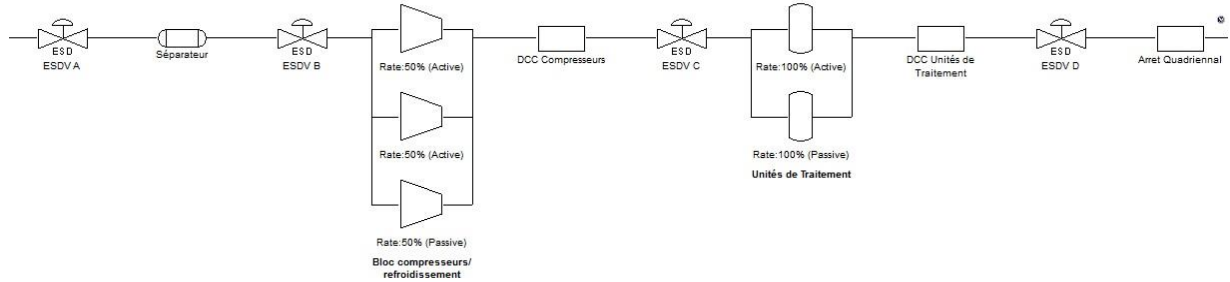

**Figure 1.** Modèle MAROS : Bloc diagramme de flux général

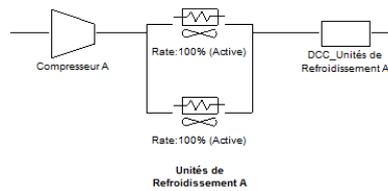

**Figure 2.** Modèle MAROS : bloc diagramme de flux des branches « compresseurs/refroidissement »

| Equipment | Description | Cap Loss Fail | Cap Loss Repair | Mode | Failure Type | Failures (per 10^t) | Repair Type | MTTR |
|---|---|---|---|---|---|---|---|---|
| Compresseur A | Critical | 100% | 100% | CRT | Exponential | 1026.76 | Exponential | 17 |
| Compresseur A | Degraded | 50% | 100% | DEG | Exponential | 1084.4 | Exponential | 10 |

**Figure 3.** Modèle MAROS : modes de défaillance d'un compresseur

Avec MAROS, la modélisation des défaillances de cause commune (DCC) est obtenue par l'ajout de blocs, placés en série des éléments concernés (cf. « DCC » sur les Figures 1 et 2), selon la même logique que les blocs diagrammes de fiabilité. Les maintenances préventives et l'arrêt général sont modélisés par des évènements déterministes (définis par des fréquences et des durées), placés en série des éléments concernés (cf. « arrêt quadriennal » sur la Figure 1, les maintenances préventives sont quant à elles des sous-blocs des blocs équipements). Les maintenances préventives peuvent également être définies comme étant « simultanées » à l'« arrêt quadriennal ». Les équipes et ressources de maintenance sont définies pour chaque équipement tel que présenté sur la Figure 4 (exemple d'un compresseur). Les caractéristiques des équipes et ressources de maintenance sont quant à elles définies par des blocs spécifiques tel que présenté sur la Figure 5.

**Figure 4.** Modèle MAROS : équipes et ressources de maintenance pour un compresseur

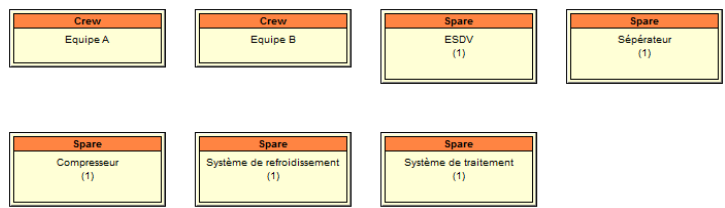

**Figure 5.** Modèle MAROS : définition des équipes et ressources de maintenance

**IV.3. Application de BStok au cas d'étude**

La Figure 6 présente le bloc diagramme de flux du modèle BStok. À noter que chaque bloc peut aussi être représenté par une image, ce qui rend alors l'affichage du modèle plus convivial. Pour chaque équipement, les modes de défaillance et maintenances préventives sont définis tel que présenté sur les Figures 7 et 8 (exemple d'un compresseur).

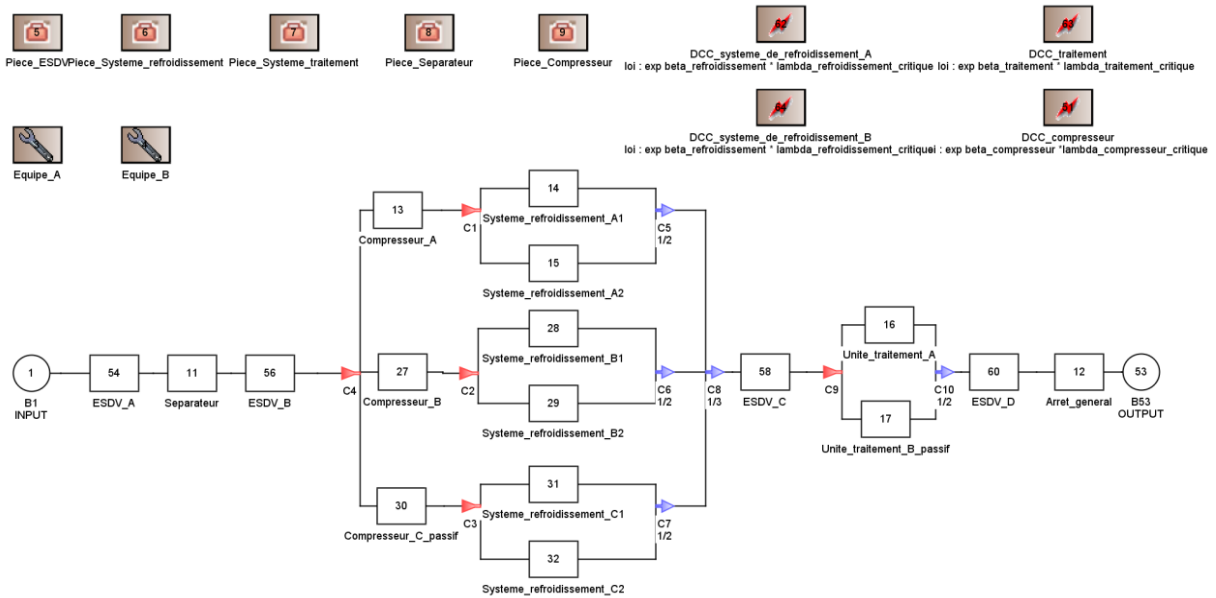

**Figure 6.** Modèle BStok : Bloc diagramme de flux général

**Figure 7.** Modèle BStok : modes de défaillance d'un compresseur

**Figure 8.** Modèle BStok : maintenances préventives d'un compresseur

Avec BStok, les défaillances de cause commune (DCC), équipes et ressources de maintenance sont définis au sein de blocs spécifiques (cf. « Pièces », « Equipe », et « DCC » sur la Figure 6). Par exemple, les caractéristiques des ressources de maintenance sont alors définies tel que présenté sur la Figure 8 (exemple avec les pièces de rechange des compresseurs). L'arrêt général est modélisé par un évènement déterministe (définis par des fréquences et des durées), placés en série du bloc diagramme (cf. « arrêt général » sur la Figure 6).

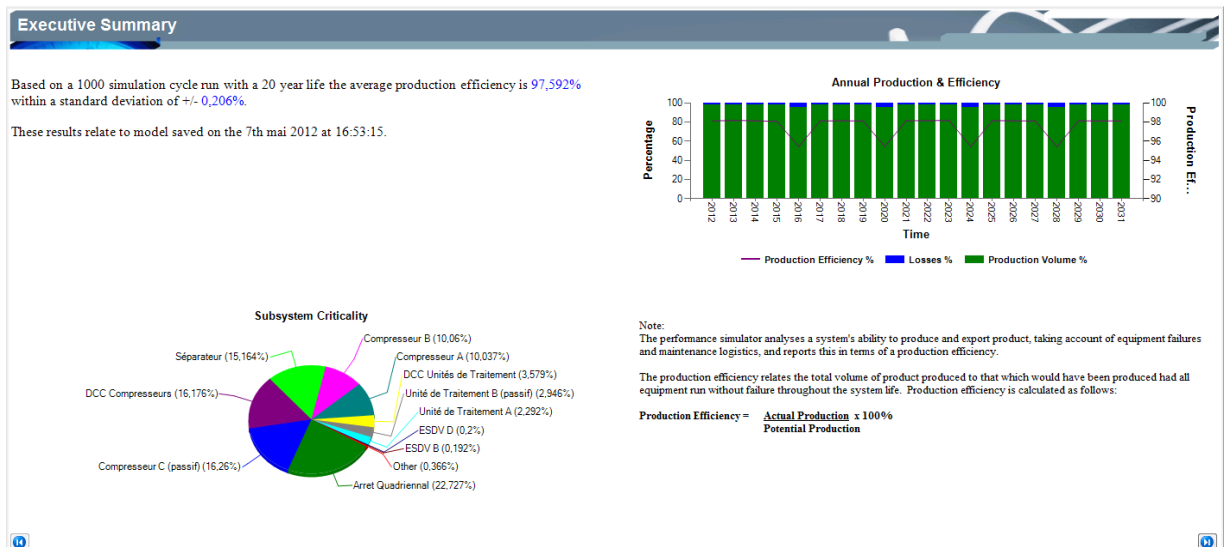

**Figure 8.** Modèle BStok : définition des pièces de rechange des compresseurs

**IV.4. Comparaison des résultats**

Pour chacun des deux outils, les analyses ont été obtenues en utilisant des séries de 1000 simulations.

MAROS intègre une interface graphique interactive, tel que présenté sur la Figure 9, qui renseigne notamment l'utilisateur sur : la disponibilité moyenne de production ; l'évolution de la disponibilité de production ; divers indicateurs de criticité (dont les définitions mathématiques ne sont malheureusement pas précisées) ; le nombre de défaillances, temps avant défaillance, et temps de réparation moyens de chaque équipement. De plus, un suivi de l'évolution de variables optionnelles (e.g. quantité d'envoi à la torche, états des équipes et des ressources de maintenance, états des « buffers ») est également proposé.

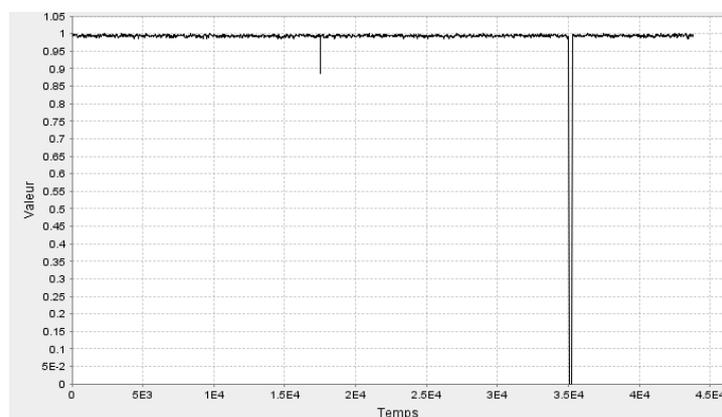

**Figure 9.** Résultats MAROS : interface graphique

Les résultats de BStok doivent quant à eux être extraits manuellement par l'intermédiaire des variables du modèle et les sept types de statistiques proposées. Par exemple, la Figure 10 représente l'évolution de la disponibilité de production sur les cinq premières années. Sur cette figure, le « pic de réduction » au temps 17 520 heures (i.e. 2 ans) correspond à la maintenance préventive des compresseurs, et le « pic de réduction » au temps 35 040 heures (i.e. 4 ans) correspond à l'arrêt général.

**Figure 10.** Résultats BStok : évolution de la disponibilité de production sur les 5 premières années

Les simulations ont permis d'obtenir une disponibilité moyenne de production (sur 20 ans) de 96,95% avec BStok (avec un écart-type de 0.31%), et de 97,59% avec MAROS (avec un écart-type de 0.21%), soit un écart absolu de 0,64%. Afin d'expliquer cette différence, des simulations complémentaires ont été effectuées. Notamment, il a été constaté que la suppression des temps de mobilisation des équipes et de réapprovisionnement des ressources de maintenance (i.e. temps réduits à zéro) a permis d'obtenir un écart absolu réduit à 0,08% pour la disponibilité moyenne de production (avec un résultat de 98,56% pour BStok, et de 98,64% pour MAROS). Ces observations suggèrent donc que ces deux outils se basent sur des modalités spécifiques de gestion de la maintenance (par exemple, le temps de mobilisation est-il compter deux fois pour la réparation, sans discontinuité, de deux équipements ? Les maintenances préventives et pertes associées sont-elles correctement reportées en cas d'indisponibilité des équipes de maintenance ?).

Une analyse de disponibilité de production peut également être utilisée pour identifier des mesures d'amélioration des performances. Pour cela, les indicateurs de performance suivants peuvent être définis (Brissaud *et al.*, 2012) :
- La criticité d'un sous-système, définie comme étant égale à 100% moins la disponibilité moyenne de production (globale) lorsque la capacité du sous-système est maintenue à zéro (et que les maintenances qui le concerne n'ont pas lieu).
- La contribution d'un sous-système à l'indisponibilité moyenne de production (ou « potentiel d'amélioration » à l'indisponibilité moyenne de production), définie comme étant égale à la disponibilité moyenne de production (globale) lorsque la capacité du sous-système est maintenue à 100% (i.e. le sous-système est « parfait », et il n'est notamment pas concerné par les maintenances) moins la disponibilité moyenne de production obtenue dans le cas « de base ».

Les résultats obtenus par BStok et MAROS pour ces indicateurs de performance sont reportés dans la Table 7. (À noter que ces indicateurs sont différents de ceux fournis dans l'interface des résultats de MAROS.)

**Table 7.** Indicateurs de performance

|  | Criticité | | Contribution à l'indisponibilité moyenne de production | |
|---|---|---|---|---|
|  | *BStok* | *MAROS* | *BStok* | *MAROS* |
| **ESDV** | 100,00% | 100,00% | 0,04% | 0,02% |
| **Compresseurs** | 100,00% | 100,00% | 1,88% | 1,20% |
| **Compresseur A** | 8,66% | 8,39% | 1,44% | 0,85% |
| **Compresseur B** | 8,66% | 8,39% | 1,43% | 0,85% |
| **Compresseur C** | 8,77% | 8,39% | 1,38% | 0,84% |
| **Unités de refroidissement** | 100,00% | 100,00% | 0,06% | 0,08% |
| **Unités de traitement** | 100,00% | 100,00% | 0,32% | 0,24% |
| **Unité de traitement A** | 4,56% | 4,01% | 0,29% | 0,15% |
| **Unité de traitement B** | 4,58% | 3,92% | 0,28% | 0,12% |
| **Séparateur** | 100,00% | 100,00% | 0,32% | 0,35% |

Pour la criticité, les écarts relatifs notables dans les résultats obtenus sont compris entre 3% et 4% pour les compresseurs, et entre 12% et 14% pour les unités de traitement. Ces écarts, relativement faibles, peuvent en partie s'expliquer par les différences de gestion de maintenance des deux outils.

Pour la contribution à l'indisponibilité moyenne de production, les écarts relatifs notables dans les résultats obtenus sont compris entre 36% et 41% pour les compresseurs, et entre 25% et 57% pour les unités de traitement. En valeurs absolues, ces écarts sont nettement plus importants pour les compresseurs, ce qui se justifie par le fait que les compresseurs ont des conditions de maintenance plus sévères (le temps de mobilisation pour l'Équipe B de maintenance dédiée aux compresseurs est de 24 heures, tandis que le temps de mobilisation pour l'Équipe A de maintenance qui concerne les unités de traitement est de 1 heure), et des taux de défaillance plus importants, ce qui accentue l'effet des différences de gestion de maintenance des deux outils. La suppression des temps de mobilisation des équipes et de réapprovisionnement des ressources de maintenance permet ainsi de réduire fortement ces écarts de résultats (par exemple, pour les compresseurs, l'écart relatif maximal pour la contribution à l'indisponibilité moyenne de production est alors réduit à 8%).

Enfin, les biais dans les résultats de criticité, et ceux de contribution à l'indisponibilité moyenne de production, pour les éléments passifs (Compresseur C et Unité de traitement B) par rapport aux éléments actifs, s'expliquent par les différences de modélisation en ce qui concerne les maintenances préventives et les défaillances à la sollicitation.

## V. Conclusions

MAROS et BStok sont deux outils qui permettent de modéliser et d'analyser les performances d'un système de production.

En ce qui concerne la modélisation, MAROS et BStok offrent tous deux un certain nombre de fonctionnalités qui peuvent être utilisées efficacement pour prendre en compte les principales caractéristiques d'un système de production. Les fonctionnalités « de base » offertes par ces deux logiciels diffèrent toutefois sur plusieurs points et, afin de prendre en compte certaines caractéristiques particulières dans la modélisation, l'un ou l'autre de ces outils peut se révéler plus approprié (par exemple, la prise en compte de plusieurs produits, d'envois à la torche, de stockage, de rampe de production, et d'évènements climatiques est plus adaptée avec MAROS, tandis que la prise en compte de dépendances entre modes de défaillance, de modes de défaillance d'équipements de sécurité, de causes communes de défaillance, et de dépendance entre taux de défaillance ou de réparation et flux de produit ou capacités des éléments est plus adaptée avec BStok). Les comparaisons présentées dans cette communication pourraient ainsi constituer quelques pistes de développement pour l'un et l'autre de ces outils.

En ce qui concerne les analyses, MAROS et BStok permettent tous deux de fournir des résultats cohérents qui offrent ainsi de précieux indicateurs pour prédire les performances de production d'un système, identifier les éléments les plus critiques et trouver des mesures d'amélioration des performances, et comparer différentes options. Certains écarts non négligeables ont cependant pu être observés entre les résultats fournis par l'un et l'autre de ces outils. Des modalités spécifiques de gestion de la maintenance ont notamment pu être mises en évidence, ce qui explique en partie les écarts de résultats observés, mais laissent place à d'autres interrogations concernant les discordances dans les hypothèses sous-jacentes de ces deux outils.

Enfin, si l'utilisation de MAROS offre plus de « convivialité », une contrepartie est l'aspect « boite noire » qui peut en résulter, notamment en ce qui concerne les résultats fournis. BStok, grâce à l'exploitation sous-jacente des réseaux de Petri, peut offrir une plus grande flexibilité de modélisation, moyennant la nécessité de disposer d'un certain degré d'expertise et d'un temps d'étude plus important. Alors que MAROS peut donc être jugé très efficace pour les cas « classiques » d'analyses de disponibilité de production dans le domaine du pétrole et du gaz, BStok, associé aux réseaux de Petri sous-jacents, est probablement plus adapté aux cas particuliers et plus « complexes ».

## VI. **Références**

Aven T., 1985, "Reliability evaluation of multistate systems with multistate components," *IEEE Transactions on Reliability*, vol. R-34(5), p. 473-479, 1985.

Aven T., 1987, "Availability evaluation of oil/gas production and transportation systems," *Reliability Engineering*, vol. 18(1), p. 35-44, 1987.

Boiteau M., Dutuit Y., Rauzi A., Signoret J.P., 2006, "The AltaRica data-flow language in use: modeling of production availability of a multistates system," *Reliability Engineering and System Safety*, vol. 91(7), p. 747-755, 2006.

Brissaud F., Varela H., Declerck B., Bouvier N., "Production availability analysis for oil and gas facilities: Concepts and procedure," *Actes des conférences PSAM11 et ESREL 2012*, 25-29 Juin 2012, Helsinki, France (accepté pour publication).

CEI, 2001, *IEC 60050-191, Vocabulaire Électrotechnique International – Chapitre 191 : Sûreté de fonctionnement et qualité de service*, 2ème édition, Genève: Commission Électrotechnique Internationale, 2001.

Chang K.P., Chang D., Zio E., 2010, "Application of Monte Carlo Simulation for the Estimation of Production Availability in Off-Shore Installations," J. Faulin, A.A. Juan, S. Martorell, J.E. Ramirez-Márquez (éditeurs), *Simulation Methods for Reliability and Availability of Complex Systems*, Chapitre 11, Londres : Springer Verlag, 2010.

DNV, 2012, http://www.dnvusa.com/services/software/products/safeti/safetipf/maros.asp

GRIF, 2012, http://grif-workshop.com/grif/bstok-module/

ISO, 2008, *ISO 20815: Petroleum, petrochemical and natural gas industries – Production assurance and reliability management*, 1ère édition, Geneva: Organisation internationale de normalisation, 2008.

NORSOK, 1998, *NORSOK Z-016: Regularity Management & Reliability Technology*, 1ère édition, Oslo: Norwegian Technology Standards Institution, 1998.

Signoret J.P., Boiteau M., Rauzy A., Thomas P., 2004, "Production availability: new tools are coming!," *Actes du congrès Lambda-Mu 14*, 12-14 Octobre 2004, Bourges, France.

Zio E., Baraldi P., Patelli E., 2006, "Assessment of the availability of an offshore installation by Monte Carlo simulation," *International Journal of Pressure Vessels and Piping*, vol. 83(4), p. 312-320, 2006.